\documentclass[runningheads,a4paper]{llncs}

\usepackage{amssymb}
\usepackage{graphicx}
\usepackage{makeidx}  
\usepackage{amsmath}
\usepackage{amsfonts}
\begin{document}
\mainmatter              
\title{A hybrid exact-ACO algorithm for the joint scheduling, power and cluster assignment in cooperative wireless networks
\thanks{This is the author's final version of the paper published in G. Di Caro, G. Theraulaz (eds.), BIONETICS 2012: Bio-Inspired Models of Network, Information, and Computing Systems. Lecture Notes of the Institute for Computer Sciences, Social Informatics and Telecommunications Engineering (LNICST, vol. 134, pp. 3-17. Springer, Heidelberg, 2014, DOI: 10.1007/978-3-319-06944-9\_1 ). The final publication is available at Springer via http://dx.doi.org/10.1007/978-3-319-06944-9\_1.
This work was partially supported by the \emph{German Federal Ministry of Education and Research} (BMBF), project \emph{ROBUKOM} \cite{BlDAHa11,KoEtAl12}, grant 03MS616E. The Author thanks Antonella Nardin for fruitful discussions.}%
}
\titlerunning{A hybrid exact-ACO algorithm for Cooperative Wireless Network Design}  
%
\author{Fabio D'Andreagiovanni\inst{1,}\inst{2}}
\authorrunning{Fabio D'Andreagiovanni}   
%
\tocauthor{Fabio D'Andreagiovanni}
\institute{Department of Optimization,\\
Zuse-Institut Berlin (ZIB), 14195 Berlin, Germany
\and
Department of Computer, Control and Management Engineering,\\
Sapienza Universit\`{a} di Roma, 00185 Roma, Italy\\
\email{d.andreagiovanni@zib.de},\\ Home page:
\texttt{http://www.dis.uniroma1.it/$\sim$fdag/}}

\maketitle              

\begin{abstract}
Base station cooperation (BSC) has recently arisen as a promising way to increase the capacity of a wireless network. Implementing BSC adds a new design dimension to the classical wireless network design problem: how to define the subset of base stations (clusters) that coordinate to serve a user. Though the problem of forming clusters has been extensively discussed from a technical point of view, there is still a lack of effective optimization models for its representation and algorithms for its solution. In this work, we make a further step towards filling such gap: 1) we generalize the classical network design problem by adding cooperation as an additional decision dimension; 2) we develop a strong formulation for the resulting problem; 3) we define a new hybrid solution algorithm that combines exact large neighborhood search and ant colony optimization. Finally, we assess the performance of our new model and algorithm on a set of realistic instances of a WiMAX network.
\end{abstract}

\section{Introduction}
Wireless communications have experienced an astonishing expansion over the last years and network operators have consequently faced the challenge of providing higher capacity with the same limited amount of radio resources. In this context, the capacity increase has been mainly pursued through interference avoidance: by limiting interference, a higher number of users can indeed be served, thus increasing network capacity.
Interference avoidance can be first obtained by a careful setting of the power emissions of the base stations constituting the network.
In recent years, it was showed that an additional reduction in interference can be obtained by allowing \emph{cooperation} among the base stations: the central feature of cooperation is that service can be provided by a number of base stations forming a \emph{cluster}, in contrast to classical systems where service is provided by a single base station. These new wireless systems are called \emph{cooperative wireless networks} (CWN).
The advantage of a CWN is clear: a user that requires a single powerful signal can be instead served by a (small) number of signals of lower power, that reduce interference towards other users. Remarkable gain in capacity are thus possible in wireless technologies supporting cooperation, such as LTE \cite{SeToBa09} and WiMAX \cite{AnGhMu07}. On the other hand, cooperation has a cost and must be kept low: cooperating base stations must indeed coordinate and need to exchange a non-negligible amount of information. For a detailed discussion of the potentialities and advantages of adopting cooperation in wireless networks, we refer the reader to the recent work \cite{GeEtAl10}.

Cooperation introduces an additional decision dimension in the problem of designing a wireless network, thus making an already complicated problem even more complicated.
In this context, Mathematical Optimization has arisen as a very effective methodology to improve the quality of design, as proven in successful industrial collaborations with major wireless providers \cite{DAMaSa12,MaRoSm06,NaAgDo06}. Optimization is furthermore able to overcome the limitations of the trial-and-error approach commonly adopted by network engineers, by pursuing a much more efficient usage of the scarce radio resources of a network.

Tough CWNs have attracted a lot of attention from a technical and theoretical point of view, there is still a lack of effective optimization models and algorithms for their design. Moreover, available models and algorithms do not treasure previous experience about how tackling specific computational issues that affect the general design problem. The main objective of our work is to make further steps towards closing such gap. Specifically, our original contributions are: 1) a new binary linear model for the design of cooperative wireless networks, that jointly optimizes activated service links, power emissions and cooperative clusters; 2) a strengthening procedure that generates a stronger optimization model, in which sources of numerical issues are completely eliminated; 3) a hybrid solution algorithm, based on the combination of a special exact large neighborhood search called RINS \cite{DaRoLP05} with ant colony optimization \cite{DoMaCo96} (note that tough generic heuristics and pure and hybrid nature-inspired metaheuristics are not a novelty in (wireless) network design - see e.g., \cite{AtEtAl10,BlEtAl11,ChKiKi08,DA11b} - to the best of our knowledge such combination has not yet been investigated). This algorithm is built to effectively exploit the precious information coming from the stronger model and can rapidly find solutions of very satisfying quality.

To our best knowledge, our model is the first that considers the joint optimization of activated service links, power emissions and cooperative clusters. Models available in literature have indeed focused attention on more specific versions of the problem. In \cite{WeEtAl11}, a set covering formulation for selecting clusters is introduced and solved by a greedy heuristic. A weakness of this work is that the feasible clusters must be explicitly listed and their number is in general exponential. In \cite{GiEtAl10}, a linear binary formulation for the problem of forming cooperative clusters and assigning users to base stations is proposed and a limited computational experience on small-sized instances is presented.
In \cite{HoEtAl12} a non-linear model for the joint optimal clustering and beamforming in CWN is investigated. Given the intrinsic hardness of the resulting non-linear problem, the Authors propose an iterative heuristic approach that makes use of mean square errors.
In contrast to our work, all these works do not consider the computational issues associated with the coefficients representing signal propagation and thus their models are weak and their solution approaches may lead to coverage errors, as explained in \cite{CaEtAl11,DAMaSa12,KeOlRa10}.

The remainder of this paper is organized as follows: in Section \ref{sec:CWND}, we introduce our new model for the design of cooperative wireless networks and discuss its strength; in Section \ref{sec:ACO} we present our hybrid metaheuristic; in Section \ref{sec:results} we present preliminary computational results over a set of realistic WiMAX instances and, finally, in Section \ref{sec:conclusions}, we derive some conclusions.

\section{Cooperative Wireless Network Design}
\label{sec:CWND}
We consider the design of a wireless network made up of a set $B$ of base stations (BSs) that provide a telecommunication service to a set $T$ of user terminals (UTs). A UT $t \in T$ is said to be \emph{covered} (or \emph{served}), if it receives the service within a minimum level of quality.
The BSs and the UTs are characterized by a number of parameters (e.g., geographical location, power emission, frequency used for transmitting). The \emph{Wireless Network Design Problem (WND)} consists in setting the values of the base station parameters, with the aim of maximizing a revenue function associated with coverage (e.g., number of covered UTs). For an exhaustive introduction to the WND, we refer the reader to \cite{DAMaSa12,DA11a,KeOlRa10}.

An optimization model for the WND typically focuses attention only on a subset of the parameters characterizing a BS. In particular, the majority of the models considers the power emission and the assignment of served UTs to BSs as the main decision variables. These are indeed two critical decisions that must be taken by a network administrator, as indicated in several real studies (e.g., \cite{AmEtAl06,BuDA12a,BuDA12b,DAMaSa12,MaRoSm06,NaAgDo06}). Other parameters that are commonly considered are the frequency and the transmission scheme used to serve a terminal (e.g., \cite{ClKoSh13,DAMa08,MoSmAl04}). Different versions of the WND are discussed in \cite{DA11a}, where a hierarchy of WND problems is also presented.

In this work, we consider a generalization of the \emph{Scheduling and Power Assignment Problem} (SPAP), a version of the WND that is known to be HP-Hard \cite{MaRoSm06}. In the SPAP, two decisions must be taken: 1) setting the power emission of each BS and 2) assigning served UTs to activated cluster of BSs (note that this corresponds to identify a subset of service links BS-UT that can be scheduled simultaneously without interference, so we use the term \emph{scheduling}).
Since we address the problem of designing a cooperative network, in our generalization we also include a third decision: for each served UT, we want to select a subset of BSs forming the cluster that serves the UT. We call the resulting generalization by the name \emph{Scheduling, Power and Cluster Assignment Problem} (SPCAP).

To model these three decisions, we introduce three types of decision variables:
\begin{enumerate}
   \item a non-negative discrete variable $p_b \in {\cal P} = \{P_1, \ldots, P_{|{\cal P}|}\}$, with
       $P_{|{\cal P}|} = P_{\max}$ and $P_l > P_{l-1} >0$, for $l = 2, \dots, |{\cal P}|$, representing the power emission of a BS $b$. Such variable can be represented as the linear combination of the power values $P_l$ multiplied by binary variables: to this end, for each $b \in B$ we introduce one binary variable $z_{bl}$ (\emph{power variable}) that is equal to $1$ if $b$ emits at power $P_l$ and $0$ otherwise. Additionally, we must express the fact that each BS may emit only at a single power level. If we denote the set of feasible power levels by $L$, the definition of the binary power variables is: $p_b = \sum_{l \in L} P_l \hspace{0.1cm} z_{bl}$, $\forall b \in B$, with $\sum_{l \in L} z_{bl} \leq 1$.
       We remark that a BS $b$ such that $\sum_{l \in L} z_{bl} = 0$ has null power (i.e., $p_b = 0$) and therefore is not deployed in the network. Our model thus also decides the localization of the BSs to be deployed, though this feature is not explicitly pointed out by defining a specific decision variable.
   \item a binary \emph{service assignment variable} $x_{t} \in \{0,1\}$, $\forall \hspace{0.1cm} t \in T$, that is equal to 1 if UT $t \in T$ is served and to 0 otherwise.
   \item a binary \emph{cluster assignment variable} $y_{tb} \in \{0,1\}$, $\forall \hspace{0.1cm} t \in T, b \in B$, that is equal to 1 if BS $b$ belongs to the cluster that serves UT $t$ and to 0 otherwise.
\end{enumerate}

\noindent
Every UT $t \in T$ picks up signals from each BS $b \in B$ in the network and the
power $P_{tb}$ that $t$ gets from $b$ is proportional to the emitted power $p_b$ by a factor $a_{tb} \in [0,1]$, i.e. $P_{tb} = a_{tb} \cdot p_b$. The factor $a_{tb}$ is an attenuation coefficient that summarizes the reduction in power experienced by a signal propagating from $b$ to $t$ \cite{NaAgDo06}. The BSs whose signals are picked up by $t$ can be subdivided into \emph{useful} and \emph{interfering}, depending on whether they contribute either to guarantee or to destroy the service.
Given a cluster of BSs $C_t \subseteq B$ that cooperate to serve $t$, the BSs in $C_t$ provide useful signals to $t$, while all the BSs that do not belong to $C_t$ are interferers. For a given service cluster $C_t$, a UT $t$ is considered served if the ratio of the sum of the useful powers to the sum of the interfering powers (\emph{signal-to-interference ratio} or \emph{SIR}) is above a threshold $\delta_t >0$, which depends on the desired quality of service \cite{NaAgDo06}:
\begin{equation}
\label{eq:firstSIRineq}
\frac{
\sum_{b \in C_t
} a_{tb}
\hspace{0.1cm} p_b
}
{
N + \sum_{b \in B \backslash C_t
} a_{tb}
\hspace{0.1cm} p_b}
\hspace{0.1cm}
\geq
\hspace{0.1cm}
\delta_t \;  .
\end{equation}

\noindent
Note that in the denominator, we highlight the presence of the system noise $N > 0$ among the interfering terms.
By simple linear algebra operations, we can reorganize the ratio \eqref{eq:firstSIRineq} into the following inequality, commonly called \emph{SIR inequality}:
\begin{equation}
\label{eq:secondSIRineq}
\sum_{b \in C_t} a_{tb}
\hspace{0.1cm} p_b
- \delta_t \sum_{b \in B \backslash C_t} a_{tb}
\hspace{0.1cm} p_b
\hspace{0.1cm}
\geq
\hspace{0.1cm}
\delta_t \hspace{0.1cm} N \;  .
\end{equation}

\noindent
We now explain how to modify the basic SIR inequality \eqref{eq:secondSIRineq}, in order to build our optimization problem.
As first step, we need to introduce the service assignment variable $x_t$ into \eqref{eq:secondSIRineq}, as we want to decide which UTs $t \in T$ are served:
\begin{equation}
\label{eq:SIRcnstr_1}
\sum_{b \in C_t} a_{tb}
\hspace{0.1cm} p_b
- \delta_t \sum_{b \in B \backslash C_t} a_{tb}
\hspace{0.1cm} p_b
+ M \hspace{0.1cm} (1 - x_t)
\hspace{0.1cm}
\geq
\hspace{0.1cm}
\delta_t \hspace{0.1cm} N \;  .
\end{equation}

\noindent
This inequality contains a large positive constant $M$ (the so-called \emph{big-M} coefficient \cite{DAMaSa12,NeWo88}), that in combination with the binary variable $x_t$ either activates or deactivates the constraint. It is straightforward to check that if $x_t = 1$, then \eqref{eq:SIRcnstr_1} reduces to \eqref{eq:secondSIRineq} and the SIR inequality \emph{must be satisfied}. If instead $x_t = 0$ and $M$ is sufficiently large (for example, $M =  + \delta_t \hspace{0.1cm} \sum_{b \in B \backslash C_t} a_{tb} \hspace{0.1cm} P_{\max} + \delta_t \hspace{0.1cm} N$), then \eqref{eq:SIRcnstr_1} \emph{becomes redundant}, as it is satisfied by any power configuration of the BSs.
\\
As second step, we need to introduce the cluster assignment variables $y_{tb}$ into \eqref{eq:SIRcnstr_1}, as we want to decide which BSs $b \in B$ form the cluster $C_t$ serving $t$:
\begin{equation}
\label{eq:SIRcnstr_2}
\sum_{b \in B} a_{tb}
\hspace{0.1cm} p_b \hspace{0.1cm} y_{tb}
\hspace{0.1cm}
- \delta_t \sum_{b \in B} a_{tb}
\hspace{0.1cm} p_b \hspace{0.1cm} (1 - y_{tb})
\hspace{0.1cm}
+ M \hspace{0.1cm} (1 - x_t)
\hspace{0.1cm}
\geq
\hspace{0.1cm}
\delta_t \hspace{0.1cm} N \;  .
\end{equation}

\noindent
In this case, we extend the summation of the inequalities to all the BSs $b \in B$ and the service cluster $C_t$ is made up of the BSs $b$ such that $y_{tb} = 1$. The BSs that are not in the cluster act as interferers and this is expressed by including the complement $1 - y_{tb}$ of the cluster variables into the interfering summation.
By $p_b = \sum_{l \in L} P_l \cdot z_{bl}$ and by simply reorganizing the terms, we finally obtain:
\begin{equation}
\label{eq:SIRcnstr_3}
\left(1 + \delta_h \right)
\sum_{b \in B} \sum_{l \in L} a_{tb}
\hspace{0.1cm} P_l \hspace{0.1cm} z_{bl} \hspace{0.1cm} y_{tb}
\hspace{0.1cm}
- \delta_t \sum_{b \in B} \sum_{l \in L} a_{tb}
\hspace{0.1cm} P_l \hspace{0.1cm} z_{bl}
\hspace{0.1cm}
+ M \hspace{0.1cm} (1 - x_t)
\hspace{0.1cm}
\geq
\hspace{0.1cm}
\delta_t \hspace{0.1cm} N \;  .
\end{equation}

\noindent
This last inequality constitutes the core of our model and is called \emph{SIR constraint}. It is actually non-linear as it includes the product $z_{bl} \hspace{0.1cm} y_{tb}$ of binary variables. However, this does not constitute an issue, as this is a special non-linearity that can be linearized in a standard way, by introducing a new binary variable $v_{tbl} = z_{bl} \hspace{0.1cm} y_{tb}$ for each triple $(t,b,l): t \in T, b \in B, l \in L$, and three new constraints: (c1) \hspace{0.05cm} $v_{tbl} \leq z_{bl}$; \hspace{0.05cm} (c2) \hspace{0.05cm} $v_{tbl} \leq y_{tb}$, \hspace{0.05cm} (c3) \hspace{0.05cm} $v_{tbl} \geq z_{bl} + y_{tb} -1$ \cite{HaRu68}.
If $z_{bl} = 0$ or $y_{tb} = 0$, then by any of (c1),(c2) we have $v_{tbl} = 0$. If instead $z_{bl} = y_{tb} = 1$, then by (c3) we have $v_{tbl} = 1$.
After these considerations, we can finally state the following pure binary \emph{linear} formulation for the SPCAP:
\begin{eqnarray}
    &&
    \max \hspace{0.1cm} \sum_{t \in T} r_{t} \hspace{0.1cm} x_{t}
    - \sum_{t \in T}  c_t \hspace{0.1cm}  \left(\sum_{b \in B} y_{tb} - x_t\right)
    \hspace{3.5cm}
    \mbox{(big-M SPCAP)}
    \nonumber
    \\
    \nonumber
    \\
    &&
    \left(1 + \delta_h \right)
    \sum_{b \in B} \sum_{l \in L} a_{tb}
    \hspace{0.1cm} P_l \hspace{0.1cm} v_{tbl}
    \hspace{0.05cm}
    - \delta_t \sum_{b \in B} \sum_{l \in L} a_{tb}
    \hspace{0.1cm} P_l \hspace{0.1cm} z_{bl}
    \hspace{0.05cm}
    + M \hspace{0.1cm} (1 - x_t)
    \geq
    \delta_t \hspace{0.1cm} N
    \nonumber
    \\
    &&
    \hspace{7.7cm}
    t \in T,
    \label{PSCAP_SIR}
    \\
    &&
    \sum_{l \in L} z_{bl} \leq 1
    \hspace{6.05cm}
    b \in B,
    \label{PSCAP_GUB}
    \\
    &&
    v_{tbl} \leq z_{bl}
    \hspace{6.3cm}
    t \in T, b \in B, l \in L,
    \label{PSCAP_lineariz1}
    \\
    &&
    v_{tbl} \leq y_{tb}
    \hspace{6.3cm}
    t \in T, b \in B, l \in L,
    \label{PSCAP_lineariz2}
    \\
    &&
    v_{tbl} \geq z_{bl} + y_{tb} -1
    \hspace{4.85cm}
    t \in T, b \in B, l \in L,
    \label{PSCAP_lineariz3}
    \\
    &&
    x_{t} \in \{0,1\}
    \hspace{6.05cm}
    t \in T,
    \label{PSCAP_varX}
    \\
    &&
    z_{bl} \in \{0,1\}
    \hspace{5.95cm}
    b \in B, l \in L,
    \label{PSCAP_varZ}
    \\
    &&
    y_{tb} \in \{0,1\}
    \hspace{5.95cm}
    t \in T, b \in B,
    \label{PSCAP_varY}
    \\
    &&
    v_{tbl} \in \{0,1\}
    \hspace{5.85cm}
    t \in T, b \in B, l \in L \; .
    \label{PSCAP_varV}
\end{eqnarray}

\noindent
The objective function maximizes the difference between the total revenue obtained by serving UTs (each served UT $t \in T$ grants a revenue $r_t > 0$) and the total cost incurred to establish cooperation between BSs (for each UT $t$, the cooperation cost $c_t > 0$ arises for each cooperating BS when at least two BSs are cooperating; no cost arises when the service is provided by a single BS). The \emph{linear} SIR constraints \eqref{PSCAP_SIR} express the coverage conditions under cooperation. The constraints \eqref{PSCAP_GUB} ensure that each BS emits at a single power level, while the constraints \eqref{PSCAP_lineariz1}-\eqref{PSCAP_lineariz3} operate the linearization of the product $z_{bl} \hspace{0.1cm} y_{tb}$. Finally, \eqref{PSCAP_varX}-\eqref{PSCAP_varV} define the decision variables of the problem.

\subsubsection{Strengthening the formulation for the PSCAP.}
The formulation (big-M SPCAP) constitutes a natural optimization model for the SPCAP. However, because of the presence of the big-M coefficients, it is known to be really weak, i.e. its linear relaxation produces very low quality bounds \cite{NeWo88}. Moreover, as the fading coefficients typically vary in a very wide range, the coefficient matrix may be very ill-conditioned and this leads to heavy numerical instability in the solution process. As a consequence, the effectiveness of standard solution algorithms provided by state-of-the-art commercial solvers, such as IBM ILOG Cplex \cite{CPLEX}, may be greatly reduced and solutions may even contain coverage errors, as pointed out in several works (e.g., \cite{CaEtAl11,DAMaSa12,DA11a,DAMaSa11,KeOlRa10,MaRoSm06}).

In order to tackle these computational issues, we extend a very effective strengthening approach that we proposed in \cite{DAMaSa12,DA11a,DAMaSa09}: the basic idea is to exploit the \emph{generalized upper bound (GUB)} constraints $\sum_{l \in L} v_{tbl} \leq 1$, implied by the other GUB constraints $\sum_{l \in L} z_{bl} \leq 1$ (as $v_{tbl} \leq z_{bl}$), to replace the SIR constraints \eqref{PSCAP_SIR} with a set of \emph{GUB cover inequalities}. For an exhaustive introduction to the cover inequalities and to their GUB version, we refer the reader to the book \cite{NeWo88} and to \cite{Wo90}. Here, we briefly recall the well-known main results about the general \emph{cover inequalities}: a \emph{knapsack constraint} $\sum_{j \in J} a_j x_j \leq b$ with $a_j, \hspace{0.05cm} b \in \mathbb{R}_+$ and $x_j \in \{0,1\}$, $\forall j \in J$, can be replaced with a (in general exponential) number of cover inequalities $\sum_{j \in C} x_j \leq |C| - 1$, where $C$ is a \emph{cover}. A cover is a subset of indices $C \subseteq J$ such that the summation of the corresponding coefficients $a_j, j \in C$ violates the knapsack constraint, i.e. $\sum_{j \in C} a_j > b$. The cover inequalities thus represent combinations of binary variables $x_j$ that cannot be activated at the same time (therefore, we can activate at most $|C| - 1$ variables in each cover $C$). The GUB cover inequalities are a stronger version of the simple cover inequalities, that can be defined when there are additional constraints of the form $\sum_{j \in K \subseteq J} x_j \leq 1$ (GUB constraints), that allows to activate at most one variable in the subset $K$.

By applying the definition of \cite{Wo90} and reasoning similarly to \cite{DAMaSa12}, we can define the general form of the GUB cover inequalities (GCIs) needed to replace the SIR constraint \eqref{PSCAP_SIR}:
\begin{eqnarray}
\label{GUB_SIR}
    x_t + \sum_{i = 1}^{|\Delta|} \sum_{l = 1}^{\lambda_i} v_{tbl}
    + \sum_{i = 1}^{|\Gamma|} \sum_{l = q_i}^{|L|} z_{bl}
    \leq |\Delta| + |\Gamma| \; ,
\end{eqnarray}
with $t \in T$, \hspace{0.05cm} $\Delta \subseteq B$, \hspace{0.05cm}
$\lambda = (\lambda_1,\ldots,\lambda_{|\Delta|}) \in L^{|\Delta|}$, \hspace{0.05cm}
$\Gamma \subseteq B \backslash \Delta$, $(q_1,\ldots,q_{|\Gamma|}) \in L^I(t,\Delta,\lambda,\Gamma)$, with $L^I(t,\Delta,\lambda,\Gamma) \subseteq L^{|\Gamma|}$ representing the subset of interfering levels of BSs in $\Gamma$ that destroy the service of $t$ provided by the cluster $\Delta$ of BSs, emitting with power levels $\lambda = (\lambda_1,\ldots,\lambda_{|\Delta|})$.
Intuitively, for fixed UT, subset of serving BSs and subset of interfering BSs, a GCI is built by fixing a power setting of the serving BSs and defining a power setting of the interfering BSs that deny the coverage of the considered UT.

If we replace the big-M SIR constraints \eqref{PSCAP_SIR} with the GCIs \eqref{GUB_SIR} in the big-M formulation (big-M SPCAP), we obtain what we call a \emph{Power-Indexed} formulation (PI-SPCAP).
The \emph{Power-Indexed} formulation (PI-SPCAP) does not contain big-M and fading coefficients and is very strong and completely stable. On the other hand, it generally contains an exponential number of constraints and must be solved through a typical \emph{branch-and-cut} approach \cite{NeWo88}: initially we define a starting formulation containing only a suitable subset of GCIs \eqref{GUB_SIR}, then we insert additional GCIs if needed, through the solution of an auxiliary separation problem (we refer the reader to \cite{DAMaSa12} for details about the separation of the GCIs of a Power-Indexed formulation).
In our case, the starting formulation contains the GCIs:
\begin{eqnarray}
\label{relaxed_GUB}
x_t + \sum_{l = 1}^{\lambda} v_{t \beta l}
+ \sum_{l = q}^{|L|} z_{bl}
\leq 2.
\end{eqnarray}

\noindent
Such GCIs are obtained by considering a relaxed version of the SIR constraints \eqref{PSCAP_SIR}, which contain only a single-server and a single-interferer (i.e., there are no summations over multiple serving and interfering BSs). This relaxation comes from the practical observation that it is common to have a server and an interferer BSs that are sensibly stronger than all the other and coverage of the user just depends on the power configuration of them \cite{DAMaSa12,DAMaSa11}. Computational experience also shows that the relaxed SIR constraints provide a good approximation of (big-M SPCAP) \cite{DAMaSa12}.
The GCIs \eqref{relaxed_GUB}
of the relaxed SIR constraints can be used to define a relaxed Power-Indexed formulation denoted by (PI$^0$-SPCAP), that constitutes a very good starting point for a branch-and-cut algorithm used to solve (PI-SPCAP), as reported in \cite{DAMaSa12}.

Following the features of the improved ANT algorithm proposed in \cite{Ma99},
in our ANT algorithm presented in the next section, we make use of two lower bounds for the SPCAP: 1) the one obtained by solving the linear relaxation of (PI-SPCAP), that we will denote by \emph{PI-bound}; 2) the one obtained by solving the linear relaxation of (big-M SPCAP), strengthened with the GCIs of (PI$^0$-SPCAP),
that we will denote by \emph{BM-bound};
We remark that \emph{PI-bound} is consistently better than \emph{BM-bound}, as it comes from a stronger formulation. However, its computation takes more time than that of \emph{BM-bound}, as it requires to separate additional GCIs (we recall that (PI-SPCAP) initially corresponds to (PI$^0$-SPCAP) and its solution generally requires to generate additional GCIs).

\section{A hybrid exact-ACO algorithm for the SPCAP}
\label{sec:ACO}
Ant Colony Optimization (ACO) is a metaheuristic approach to combinatorial optimization problems that was originally proposed by Dorigo and colleagues, in a series of works from the 1990s (e.g., \cite{DoMaCo96}), and it was later extended to integer and continuous optimization problems (e.g., \cite{DoDCGa99}). For an overview of the theory and applications of ACO, we refer the reader to \cite{BlEtAl11,DoBl05,DoDCGa99}.
It is now common knowledge that the algorithm draws its inspiration from the foraging behaviour of ants. The basic idea of an ACO algorithm is to define a loop where a number of feasible solutions are iteratively built in parallel, exploiting the information about the quality of solutions built in previous executions of the loop. The general structure of an algorithm can be depicted as follows:
\begin{description}
  \item UNTIL an arrest condition is not satisfied DO
      \hspace{2.0cm} (Gen-ACO)
      \begin{enumerate}
      \item Ant-based solution construction
      \item Daemon actions
      \item Pheromone trail update
    \end{enumerate}
\end{description}
We now describe the details of each phase presented above for our hybrid exact-ACO algorithm for the SPCAP. Our approach is hybrid since the canonical ACO step 1 is followed by a daemon action phase, where we exactly explore a large neighborhood of the generated feasible solutions, by exactly solving a Mixed-Integer Program, as explained in Subsection \ref{sec:RINS}.

\subsection{Ant-based solution construction}
In step 1 of the loop, $m \in \mathbb{Z}_+$ computational agents called \emph{ants} are defined and each ant iteratively builds a feasible solution for the optimization problem. At every iteration, the ant is in a so-called \emph{state}, associated with a \emph{partial solution} to the problem, and can complete the solution by selecting a \emph{move} among a set of feasible ones. The move is probabilistically chosen on the basis of its associated pheromone trail values. For a more detailed description of the elements and actions of step 1, we refer the reader to the paper by Maniezzo \cite{Ma99}. The paper proposes an improved ANT algorithm (ANTS), which we take as reference for our work. We were particularly attracted by the improvements proposed in ANTS, as they are based on the attempt of better exploiting the precious information that comes from a \emph{strong} Linear Programming formulations of the original discrete optimization problem. Moreover, ANTS also uses a reduced number of parameters and adopts mathematical operations of higher computational efficiency (e.g., multiplications instead of powers with real exponents).

Before describing the structure and the behaviour of our ants, we make some preliminary considerations. Our formulation for the SPCAP is based on four types of variables: 1) power variables $z_{bl}$; 2) cluster variables $y_{tb}$; 3) service variables $x_{t}$; 4) linearization variables $v_{tbl}$. Once that the power variables and the cluster variables are fixed: i) the value of the linearization variables is immediately determined, because of constraints \eqref{PSCAP_lineariz1}-\eqref{PSCAP_lineariz3}, and ii) the objective function can be easily computed, as the served UTs can be identified by simply checking which SIR inequalities \eqref{PSCAP_SIR} are satisfied. As a consequence, in the ant-construction phase we can limit our attention to the cluster and service variables and we introduce the following two concepts of power and cluster states.
\begin{definition}
  Power state (PS): let $L_0 = L \cup \{0\}$ be the set of power levels plus the null power level. A \emph{power state} represents the activation of a subset of BSs on some power level $l \in L_0$ and excludes that the same BS is activated on two power levels. Formally: $PS \subseteq B \times L_0: \not \exists (b_1,l_1), (b_2,l_2) \in PS: b_1 = b_2$.
\end{definition}

\noindent
We say that a power state PS is \emph{complete} when it specifies the power configuration of every BS in B (thus $|PS| = |B|$). Otherwise the PS is called \emph{partial} and we have $|PS| < |B|$. Furthermore, for a given power state PS, we denote by $B(PS)$ the subset of BSs whose power is fixed in PS (we call such BSs \emph{configured}), i.e. $B(PS) = \{b \in B: \exists (b,l) \in PS\}$.
\begin{definition}
  Cluster state (CS): let $\bar{B} = B$
   and let $T \times \bar{B}$ be the set of couples $(t,b)$ representing the \emph{non-}assignment of BS $b$ to the cluster serving UT $t$. A \emph{cluster state} represents the assignment or non-assignment of a subset of BSs to the clusters serving a subset of UTs and excludes that the same BS is at the same time assigned and non-assigned to the cluster of a UT. Formally:
   $CS \subseteq T \times B \hspace{0.1cm} \cup \hspace{0.1cm} T \times \bar{B}:
   \hspace{0.1cm}
   \not \exists (t_1,b_1), (t_2,b_2) \in CS:
   \hspace{0.1cm}
   b_1 \in B, b_2 \in \bar{B}
   \mbox{ and } t_1 = t_2,  b_1 = b_2 = b$.
\end{definition}

\noindent
We say that a cluster state CS is \emph{complete} when it specifies the cluster assignment or non-assignment of every BS in $B$ to every UT in $T$ (thus $|CS| = |T||B|$). Otherwise the CS is called \emph{partial} and we have $|CS| < |T||B|$).
Moreover, for a given cluster state $CS$ and UT $t$, we denote by $B_t(CS)$ the subset of BSs that are either assigned or not assigned to the service cluster of $t$ in $CS$, i.e. $B_t(CS) = \{b \in B \cup \bar{B} : \exists (t,b) \in CS\}$.

In our ANT algorithm, we decided to first establish the value of the power variables and then the cluster variables. So an ant first passes through a sequence of partial power states, till a complete one is reached (power construction phase). Then it passes through a sequence of partial cluster states, till a complete one is reached (cluster construction phase). More formally, in the power phase, an ant moves from a partial power state $PS_i$ to a partial power state $PS_j$ such that:
$$
PS_j = PS_i \cup \{(b,l)\} \hspace{0.1cm} \mbox{ with } (b,l) \in B \times L_0: b \not \in B(PS_i) \;  .
$$
Note that by the definition of power state, the added couple $(b,l)$ may not contain a BS whose power is already fixed in a previous power state.

In the cluster phase, an ant moves from a partial cluster state $CS_i$ to a partial cluster state $CS_j$ such that:
$$
CS_j = CS_i \cup \{(t,b)\} \hspace{0.1cm} \mbox{ with } \{(t,b)\} \in T \times B \cup T \times \bar{B} \;  .
$$
Note that by the definition of cluster state, the added couple $(t,b)$ may not contain a BS that has been already either assigned or not assigned to the same UT.
Moreover, note that the definitions of power and cluster state can be immediately traced back to a sequence of fixing of the decision variables, thus relaxing the concept of state and reducing a move to the fixing of a decision variable $j$ after the fixing of another decision variable $i$, as it is done in \cite{Ma99}.

Every move adds a single new element to the partial solution. Once that the construction phases terminate, the value of the decision variables $(z,y)$ is fully established and, as previously noted, we can immediately derive the value of the other variables $(x,v)$, thus obtaining a complete feasible solution $(x,z,y,v)$ for the SPCAP.

The probability that an ant $k$ moves from a power (cluster) state $i$ to a more complete power (cluster) state $j$, chosen among a set $F$ of feasible power (cluster) states, is defined by the improved formula of \cite{Ma99}:
\begin{equation*}
p_{ij}^{k} = \frac{\alpha \hspace{0.1cm} \tau_{ij} + (1-\alpha) \hspace{0.1cm} \eta_{ij}}
                {\sum_{f \in F} \alpha \hspace{0.1cm} \tau_{if} + (1-\alpha) \hspace{0.1cm} \eta_{if}} \; ,
\end{equation*}

\noindent
where $\alpha \in [0,1]$ is the parameter establishing the relative importance of trail and attractiveness. Of course, the probability of infeasible moves is set to zero. As discussed in \cite{Ma99}, the trail values $\tau_{ij}$  and the attractiveness values $\eta_{ij}$ should be provided by suitable lower bounds of the considered optimization problem. In our particular case:
1) $\tau_{ij}$ is derived from the values of the variables in the solution associated with the  bound \emph{PI-bound}, provided by the linear relaxation of the Power-Indexed formulation (PI-SPCAP) (see the next subsection for the specific setting of $\tau_{ij}$);
2) $\eta_{ij}$ is equal to the optimal solution of the linear relaxation
of the big-M formulation (big-M SPCAP) strengthened with the GCIs of (PI$^0$-SPCAP) and including additional constraints to fix the value of the variables fixed in the considered state $j$.
We denote the latter bound by \emph{strongBM-bound} and we recall that this can be quickly computed and its computation becomes faster and faster as we move towards a complete state (the number of fixed variables indeed increases move after move).

As previously explained, once that an ant has finished its construction, we have a vector $(z,y)$ that can be used to derive the value of the other variables $(x,v)$ and define a complete feasible solution $(x,z,y,v)$ for the SPCAP.

\subsection{Daemon actions: Relaxation Induced Neighborhood Search}
\label{sec:RINS}
We refine the quality of the feasible solutions found through the ant-construction phase by an \emph{exact local search} in a \emph{large neighborhood}, made for each feasible solution generated by the ants. Specifically, we adopt a modified \emph{relaxation induced neighborhood search} (RINS) (see \cite{DaRoLP05} for a detailed discussion of the method). The main steps of RINS are 1) defining a neighborhood by exploiting information about some continuous relaxation of the discrete optimization problem, and 2) exploring the neighborhood  through a (Mixed) Integer Programming problem, that is optimally solved through an effective commercial solver.

Let $S^{\small \mbox{ANT}}$ be a feasible solution to the SPCAP built by an ant and let $S^{\small \mbox{PI}}$ be the optimal solution to the linear relaxation of (PI-SPCAP). Additionally, let $S^{\small \mbox{ANT}}_j, S^{\small \mbox{PI}}_j$ denote the $j$-th component of the vectors. Our modified RINS \emph{(mod-RINS)} solves a sub-problem of the big-M formulation (big-M SPCAP) strengthened with the GCIs of (PI$^0$-SPCAP) where:
\begin{enumerate}
    \item we fix the variables whose value in $S^{\small \mbox{ANT}}$ and $S^{\small \mbox{PI}}$ differs of at most $\epsilon > 0$, (i.e., $S_j = 0$ if $S^{\small \mbox{ANT}}_j = 0$ $\cap$ $S^{\small \mbox{PI}} \leq \epsilon$, $S_j = 1$ if $S^{\small \mbox{ANT}}_j = 1$ $\cap$ $S^{\small \mbox{PI}} \geq 1 - \epsilon$;
      \item set an objective cutoff based on the value of $S^{\small \mbox{ANT}}$;
      \item impose a solution time limit of $T$;
\end{enumerate}

\noindent
The time limit is set as the problem may be in general difficult to solve, so the exploration of the feasible set may need to be truncated. Note that in point 1 we generalize the fixing rule of RINS, in which $\epsilon = 0$.

\subsection{Pheromone trail update}
At the end of each construction phase $t$ of the ants, the pheromone trails $\tau_{ij}(t-1)$ are updated according to the following improved formula (see \cite{Ma99} for a detailed discussion of its elements):
\begin{equation}
\label{pheroFormula}
\tau_{ij}(t) \hspace{0.05cm} = \hspace{0.05cm} \tau_{ij}(t-1)
\hspace{0.05cm} + \hspace{0.05cm} \sum_{k=1}^{m} \tau_{ij}^k
\hspace{0.5cm} \mbox{ with } \hspace{0.1cm} \tau_{ij}^k =
\hspace{0.1cm}
\tau_{ij}(0)
\cdot \left(
1 - \frac{z_{curr}^k - LB}{\bar{z} - LB}
\right) ,
\end{equation}

\noindent
where, to set the values $\tau_{ij}(0)$ and $LB$, we solve the linear relaxation of (PI-SPCAP) and then we set $\tau_{ij}(0)$ equal to the values of the corresponding optimal decision variables and $LB$ equal to the optimal value of the relaxation. Additionally, $z_{curr}^k$ is the value of the solution built by ant $k$ and $\bar{z}$ is the moving average of the values of the last $\psi$ feasible solutions built. As noticed in \cite{Ma99}, formula \eqref{pheroFormula} substitutes a very sensible parameter, the pheromone evaporation factor, with the moving average $\psi$ whose setting is much less critical.

\bigskip

\noindent
The overall structure of our original hybrid exact-ACO algorithm is presented in Algorithm 1. The algorithm includes an outer loop repeated $r$ times. At each execution of the loop, an inner loop defines $m$ ants to build the solutions. Once that an ant finishes to build its solution, mod-RINS is applied in an attempt at finding an improvement. Pheromone trail updates are done at the end of each execution of the inner loop.

\bigskip

\noindent
\textbf{Algorithm 1. Hybrid exact-ACO for the SPCAP.}
\begin{enumerate}
  \item Compute the linear relaxation of (PI-SPCAP) and use it to initialize
        the values $\tau_{ij}(0)$.
  \item FOR $t := 1$ TO $r$ DO
      \begin{enumerate}
        \item FOR $\mu := 1$ TO $m$ DO
          \begin{enumerate}
            \item build a complete power state;
            \item build a complete cluster state;
            \item derive a complete feasible solution to the SPCAP;
            \item apply mod-RINS to the feasible solution.
          \end{enumerate}
        END FOR
        \item Update $\tau_{ij}(t)$ according to (\ref{pheroFormula}).
      \end{enumerate}
      END FOR
\end{enumerate}

\section{Computational experiments}
\label{sec:results}
We tested the performance of our hybrid algorithm on a set of 15 realistic instances of increasing size, defined in collaboration with the Technical Strategy \& Innovations Unit of British Telecom Italia (BT Italia). The experiments were made on a machine with a 1.80 GHz Intel Core 2 Duo processor and 2 GB of RAM and using the commercial solver IBM ILOG Cplex 11.1. All the instances refer to a WiMAX Network \cite{AnGhMu07} and lead to the definition of very large and hard to solve (big-M SPCAP) formulations. Even when strengthened with the GCIs of (PI$^0$-SPCAP), (big-M SPCAP) continues to constitute a very hard problem and the simple identification of feasible solutions may be a hard task even for Cplex. In particular, for most instances it was not possible to find feasible solutions within one hour of computations and, when solutions were found, they were anyway of low value (up to the 35\% of covered UTs). Our heuristic algorithm was instead able to find good quality solutions.

After a series of preliminary tests, we found that a good setting of the parameters of the heuristic is: $\alpha = 0.5$ (balance between attractiveness and trail level), $m = |B|/2$ (number of ants equal to half the number of BSs), $\psi = m = |B|/2$ (width of the moving average equal to the number of ants), $\epsilon = 0.01$ (tolerance of fixing in mod-RINS), $T = 10$ minutes (time limit in mod-RINS). Moreover, the construction loop was executed 50 times. In Table \ref{table:results}, for each instance we report its ID and size and the number $|T^*|$ of covered UTs in the best solution found by mod-RINS (showing also the percentage coverage $Cov\%$) and in the corresponding ant solution. We also report the maximum size of a cluster in the best solution. The solutions found by the hybrid algorithm have a much higher value than those found by Cplex directly applied to (big-M SPCAP) and guarantee a good level of coverage ranging from 45 to 80\%. Moreover, we note that the execution of mod-RINS is able to increase the value of the ant solution from 5 to 13\%. Finally, it is interesting to note that the size of the clusters keeps in general low, presenting a maximum dimension of 5. We consider the overall performance highly satisfying, considering that our real aim is to use the solutions generated by the algorithm to favour a warm start in an exact cutting-plane algorithm applied to the Power-Indexed formulation (PI-SPCAP). Moreover, we are confident that refinements of the components of the heuristic and further tuning of the procedure can lead to the generation of solutions of higher quality.
\begin{table}
\label{table:results}
\caption{Experimental results}
\small
\begin{center}
\begin{tabular}{@{\quad} l @{\quad} c @{\quad} r @{\quad} c @{\quad} c @{\quad} c @{\quad} c @{\quad}}
\hline
ID & $|T|$ & $|B|$ & $|T^*|$ (ACO) & $|T^*|$ (ACO+RINS) & Cov\% & Max size cluster
\\ [2pt]
\hline
I1 & 100 &   9 &  55 &  60 & 0.60 & 3
\\
I2 & 100 &  12 &  62 &  67 & 0.67 & 2
\\
I3 & 121 &   9 &  52 &  55 & 0.45 & 2
\\
I4 & 121 &  15 &  72 &  80 & 0.66 & 4
\\
I5 & 150 &  12 &  94 & 106 & 0.71 & 3
\\
I6 & 150 &  15 &  96 & 106 & 0.71 & 3
\\
I7 & 150 &  18 & 103 & 112 & 0.75 & 3
\\
I8 & 169 &  12 &  80 &  87 & 0.51 & 2
\\
I9 & 169 &  15 & 103 & 116 & 0.69 & 4
\\
I10 & 169 & 18 & 121 & 133 & 0.79 & 2
\\
I11 & 196 & 15 & 136 & 144 & 0.73 & 3
\\
I12 & 196 & 21 & 140 & 156 & 0.80 & 5
\\
I13 & 225 &  9 & 125 & 142 & 0.63 & 3
\\
I14 & 225 & 15 & 142 & 149 & 0.66 & 3
\\
I15 & 225 & 18 & 152 & 163 & 0.72 & 4
\\
\hline
\end{tabular}
\end{center}
\end{table}

\section{Conclusions}
\label{sec:conclusions}
Cooperative wireless networks have recently attracted a lot of attention, since cooperation among base stations may lead to remarkable increases in the capacity of a network and enhance the service experience of the users. Though base station cooperation has been extensively discussed from a theoretical and technical point of view, there is still a lack of effective optimization models and algorithms for its evaluation and implementation.
To make a further step towards filling such gap, in this work we have presented a new model and solution algorithm for the problem of designing a cooperative wireless network. In particular, we have generalized the classical model for wireless network design, in order to include cluster definition and assignment. We have then showed how to strengthen the model, through the use of a special class of valid inequalities, the GUB cover inequalities, that eliminate all the sources of numerical problems.
Finally, we have defined a hybrid heuristic based on the combination of ant colony optimization and relaxation induced neighborhood search, that exploits the important information provided by the relaxation of a strong formulation. Computational experiments on a set of realistic instances showed that our heuristic can find solutions of good quality, which could be used for a warm start in an exact branch-and-cut algorithm.
Future work will consist in refining the components of the heuristic (for example, by better integrating the power and cluster state moves) and in integrating the heuristic with a branch-and-cut algorithm, in order to find solutions of higher value and whose quality is precisely assessed.

%
%

%
\end{document}